\documentclass[12pt]{article}

\usepackage{amsmath}
\usepackage{amssymb}
\usepackage{graphicx}
\usepackage[left=3cm,right=1.5cm,top=2cm,bottom=2cm,bindingoffset=0cm]{geometry}
\usepackage{bbm}

\begin{document}

\begin{center}

Igor G. Pospelov, Stanislav A. Radionov\footnote{National Research University Higher School of Economics. Research group on macro-structural modeling of Russian economy. Intern Researcher; E-mail: saradionov@edu.hse.ru}

\textbf{TELEGRAPH PROCESS IN THE BOUNDED DOMAIN WITH ABSORBING LOWER BOUNDARY AND REFLECTING WITH DELAY UPPER BOUNDARY}\footnote{The reported study was supported by Russian Scientific Fund, research project No 14-11-00432.}
\end{center}

\begin{abstract}
We study the asymmetric one-dimensional telegraph process in the bounded domain. Lower boundary is absorbing and upper boundary is reflecting with delay. Point stays in the upper boundary until switch of regime occurs. We obtain the distribution of this process in terms of Laplace trasforms.
\end{abstract}



\section{Introduction}

Telegraph process was introduced in \cite{Goldstein1951} and \cite{Kac1974}. It was shown in these works that telegraph process is closely related to the well-known PDE called telegrapher's equation. Telegraph process has been applied to the number of problems in physics (\cite{Mugnai1992}, \cite{Joseph1989}, \cite{Ishimaru1989}) and economics (\cite{DiMasi1995}, \cite{Bondarenko2000}, \cite{Crescenzo2002}, \cite{Ratanov2007}, \cite{RatanovMelnikov2008}, \cite{Ratanov2010}, \cite{Lopez2012b}). A growing body of literature is dedicated to the properties of telegraph process: \cite{Orsingher1990} and \cite{Foong1994} derive the distribution of the symmetric telegraph process, its maximum and first passage time for unbounded domain, \cite{Masoliver1993} solve telegrapher's equation with reflecting and partially reflecting boundaries, \cite{Kaya2000} solves telegrapher's equation by the Adomian decomposition method, \cite{Beghin2001} and \cite{LopezRatanov2014} analyze asymmetric telegraph process on unbounded domain, \cite{Chen2008} and \cite{Huang2009} solve time-fractional telegrapher's equation including the case of bounded domain with different types of boundary conditions, \cite{Ostapenko2012} solves telegrapher's equation in the domain with variable borders.

We are interested in the dynamics of telegraph process in the bounded domain with absorbing lower boundary and reflecting  upper boundary. Unlike \cite{Masoliver1993}, we assume that reflection is not instantaneous --- point stays in the upper boundary until switch to the regime 0 occurs. That type of boundary conditions was analyzed in \cite{Balakrishnan1988} and the dynamics for the case of two boundaries of this type was described in \cite{Pogorui2010}. However, the case of mixed boundaries, to the best our knowledge, was not yet analyzed in the literature. We derive the system of two partial differential equations with boundary conditions and solve them by the method of Laplace transform. The solution is given in the terms of inverse Laplace transform of rather cumbersome functions.

\section{Statement and solution of the problem}

We consider telegraph process $A(t)$ in the bounded domain $[0,B]$ for some $B > 0$. There are two regimes $s(t) \in \{0,1\}$ defined by velocities $\mu_0 < 0$ for the regime $s=0$ and $\mu_1 > 0$ for the regime $s=1$. The rates of occurences of velocity switches are $\Lambda_0 > 0$ (for the switch from regime 0 to 1) and $\Lambda_1 > 0$ (for the opposite switch). The process starts at some point $A(0) \in [0,B]$ in some regime $s(0) \in \{0,1\}$. When $A(t)$ becomes negative for the first time, process stops. When $A(t)$ hits $B$, it stays in $B$ until the switch to the regime 0 accurs. In order to describe the evolution of the process, we introduce functions

\begin{equation}
F_s (t,A) = P(A(t) \geq A, s(t) = s).
\end{equation}

\noindent The definition of the process leads to the following description of dynamics

\begin{equation}\label{F}
F_0 \left( t+ \Delta t,A \right) = \left(1 -\Lambda_{{s}} \Delta t \right) { F_s} \left( t,A -\mu_{{s}}{ \Delta t} \right) +\Lambda_{{1-s}} \Delta t F_{1-s} \left( t,A \right), s \in \{0,1\}
\end{equation}

\noindent for small $\Delta t$, which leads to the system of partial differential equations

\begin{equation}\label{Fsyst}
\left\{ \begin{array}{l}
{\frac {\partial }{\partial t}}{F_0} \left( t,A \right) =- \mu_{{0}} 
{\frac {\partial }{\partial A}}{F_0} \left( t,A \right)   
-\Lambda_{{0}}{ F_0} \left( t,A \right) +\Lambda_{{1}}{ 
F_1} \left( t,A \right) ,\\ {\frac {\partial }{\partial t}}{ F_1}
 \left( t,A \right) =- \mu_{
{1}}  {\frac {\partial }{\partial A}}  {
 F_1} \left( t,A \right)  -\Lambda_{{1}}{ F_1} \left( t,A \right) +\Lambda_{{0}}{ F_0}
 \left( t,A \right).\\
 \end{array}
  \right.
\end{equation}

In order to derive boundary conditions, substitute $A=0$ and $A=B$ to \eqref{F} and note that $F_1 (t, -\mu_1 \Delta t) = F_1 (t,0)$ and $F0(t, B-\mu_0 \Delta t) = 0$. After some simple calculations we get

\begin{equation}\label{Boundary}
\frac{\partial}{\partial A} F_1 (t,0) = 0, F_0 (t,B) = 0.
\end{equation}

\noindent Our goal is to solve the system \eqref{Fsyst} with boundary conditions \eqref{Boundary}. To do so, we multiply both equations by $\rm e^{-\xi A}$ for some arbitrary $\xi$ and integrate them with respect to $A$:

\begin{equation}
\left\{ \begin{array}{l}
{\frac {\partial }{\partial t}}{ L_0} \left( t,\xi \right) =-\mu_{{0
}}\int_{0}^{B}\!{{\rm e}^{-\xi A}}{\frac {\partial }{\partial A}}{
 F_0} \left( t,A \right) \,{\rm d}A-\Lambda_{{0}}{ L_0} \left( t,
\xi \right) +\Lambda_{{1}}{ L_1} \left( t,\xi \right) , \\ {\frac {
\partial }{\partial t}}{ L_1} \left( t,\xi \right) =-\mu_{{1}}\int_{0
}^{B}\!{{\rm e}^{-\xi A}}{\frac {\partial }{\partial A}}{ F_1}
 \left( t,A \right) \,{\rm d}A-\Lambda_{{1}}{ L_1} \left( t,\xi
 \right) +\Lambda_{{0}}{ L_0} \left( t,\xi \right),
  \end{array}
  \right.
\end{equation}

\noindent where

\begin{equation}\label{Ldef}
{ L_s} \left( t,\xi \right) =\int_{0}^{B}\!{{\rm e}^{-\xi\,A}}{ 
F_s} \left( t,A \right) \,{\rm d}A, s \in \{0,1\}.
\end{equation}

\noindent Integrating by parts and using the upper boundary condition, we get

\begin{equation}\label{Lsyst}
\left\{ \begin{array}{l}
{\frac {\partial }{\partial t}}{ L_0} \left( t,\xi \right) =\mu_{{0
}}\phi \left( t \right) -\mu_{{0}}\xi\,{ L_0} \left( t,\xi \right) -
\Lambda_{{0}}{ L_0} \left( t,\xi \right) +\Lambda_{{1}}{ L_1}
 \left( t,\xi \right) ,\\ {\frac {\partial }{\partial t}}{ L_1} \left( 
t,\xi \right) =-\mu_{{1}}{{\rm e}^{-\xi\,B}}\omega \left( t \right) +\mu
_{{1}}\psi \left( t \right) -\mu_{{1}}\xi\,{ L_1} \left( t,\xi
 \right) -\Lambda_{{1}}{ L_1} \left( t,\xi \right) +\Lambda_{{0}}{
 L_0} \left( t,\xi \right),
   \end{array}
  \right.
\end{equation}

\noindent where

\begin{equation}
{ F_0} \left( t,0 \right) =\phi \left( t \right) ,{ F_1} \left( t
,B \right) =\omega \left( t \right) ,{ F_1} \left( t,0 \right) =\psi
 \left( t \right).
\end{equation}

\noindent Now taking Laplace transform of \eqref{Lsyst} with respect to $t$, we get

\begin{equation}\label{Lapl}
\begin{split}
& {\tilde{L}_0} \left( p,\xi \right) =-{
\frac { \tilde{\omega} \left( p \right) \Lambda
_{{1}}\mu_{{1}}{{\rm e}^{-\xi\,B}}}{ \left( p-n \right)  \left( p-m
 \right) }}+{\frac {\Lambda_{{1}}\mu_{{1}} \tilde{\psi}
 \left( p \right) }{ \left( p-n \right)  \left( p-m
 \right) }}+{\frac { \left( \mu_{{1}}\xi+p+\Lambda_{{1}} \right) \mu_{
{0}} \tilde{\phi} \left( p \right) }{ \left( 
p-n \right)  \left( p-m \right) }}+ \\ & {\frac {\Lambda_{{1}}{ L_1}
 \left( 0,\xi \right) }{ \left( p-n \right)  \left( p-m \right) }}+{
\frac { \left( \mu_{{1}}\xi+p+\Lambda_{{1}} \right) { L_0} \left( 0,
\xi \right) }{ \left( p-n \right)  \left( p-m \right) }},\\ & 
{ \tilde{L}_1} \left( p,\xi \right) =-{\frac {\tilde{\omega} \left( p \right) \mu_{{1}} \left( 
\mu_{{0}}\xi+p+\Lambda_{{0}} \right) {{\rm e}^{-\xi\,B}}}{ \left( p-n
 \right)  \left( p-m \right) }}+{\frac {\mu_{{1}} \left( \mu_{{0}}\xi+
p+\Lambda_{{0}} \right) \tilde{\psi} \left( p \right)  }{ \left( p-n \right)  \left( p-m \right) }}+{\frac {\Lambda
_{{0}}\mu_{{0}} \tilde{\phi} \left( p \right) }{ \left( p-n \right)  \left( p-m \right) }}+ \\ & {\frac { \left( 
\mu_{{0}}\xi+p+\Lambda_{{0}} \right) { L_1} \left( 0,\xi \right) }{
 \left( p-n \right)  \left( p-m \right) }}+{\frac {\Lambda_{{0}}{ 
L_0} \left( 0,\xi \right) }{ \left( p-n \right)  \left( p-m \right) }},
\end{split}
\end{equation}

\noindent where $p$ is the parameter of Laplace transform, tilde means Laplace transform of a function and

\begin{equation}\label{mn}
\begin{split}
& n=\frac{-\mu_{{0}}\xi-\mu_{{1}}\xi-\sqrt{q}-\Lambda_{{0}}-\Lambda_{{1}}}{2},m=\frac{-\mu_{{0
}}\xi-\mu_{{1}}\xi+\sqrt{q}-\Lambda_{{0}}-\Lambda_{{1}}}{2}, \\ & q={\xi}^{2}
\mu_0^{2}-2\,{\xi}^{2}\mu_{{0}}\mu_{{1}}+{\xi}^{2} \mu_1^{2}
+2\,\xi\,\Lambda_{{0}}\mu_{{0}}-2\,\xi\,\Lambda_{{0}}\mu_{{1}}-2\,\xi
\,\Lambda_{{1}}\mu_{{0}}+2\,\xi\,\Lambda_{{1}}\mu_{{1}}+ \\ & \Lambda_0^{2}+2\,\Lambda_{{0}}\Lambda_{{1}}+\Lambda_1^{2}.
\end{split}
\end{equation}

We now prove the simple

\textbf{Proposition 1}. 1. $m$ is positive for sufficiently big absolute values of $\xi$ and any values of the parameters of the model.

2. $n$ is negative for any values $\xi$ and any values of the parameters of the model.

\textbf{Proof}. 1. Expanding $m$ in a Taylor series in the neighborhood of $\xi = +\infty$, we get

\begin{displaymath}
m \approx -\mu_{{0}}\xi-\Lambda_{{0}}-{\frac {\Lambda_{{0}}\Lambda_{{1}}}{\xi\,
 \left( \mu_{{0}}-\mu_{{1}} \right) }},
\end{displaymath}

\noindent which is positive for sufficiently big positive values of $\xi$. Expanding $m$ in a Taylor series in the neighborhood of $\xi = -\infty$, we get

\begin{displaymath}
m \approx -\xi\,\mu_{{1}}-\Lambda_{{1}}+{\frac {\Lambda_{{0}}\Lambda_{{1}}}{\xi
\, \left( \mu_{{0}}-\mu_{{1}} \right) }},
\end{displaymath}

\noindent which is again positive for sufficiently big negative values of $\xi$.

2. Inequality $n<0$ may be rewritten as

\begin{equation}\label{neq0}
-\mu_{{0}}\xi-\xi\,\mu_{{1}}-\Lambda_{{0}}-\Lambda_{{1}}<\sqrt {q}.
\end{equation}

\noindent If the expression on the left side is negative, inequality is proven. Assume it is positive:
\begin{equation}\label{neq1}
\mu_{{0}}\xi+\xi\,\mu_{{1}}+\Lambda_{{0}}+\Lambda_{{1}}<0.
\end{equation}

\noindent Taking squares of both sides of \eqref{neq0}, we get 

\begin{equation}\label{neq2}
\xi (\xi \mu_{{0}}\mu_{{1}}+\Lambda_{{0}}\mu_{{1}}+\Lambda_{{1}}\mu_{{0}}) < 0.
\end{equation}

\noindent Assume $\mu_0+\mu_1 > 0$. Then \eqref{neq1} leads to 

\begin{displaymath}
\xi<\frac{-\Lambda_0-\Lambda_1}{\mu_0+\mu_1}<0.
\end{displaymath}

\noindent Hence, first multiplier in \eqref{neq2} is negative. Consider the second one:

\begin{displaymath}
\xi \mu_{{0}}\mu_{{1}}+\Lambda_{{0}}\mu_{{1}}+\Lambda_{{1}}\mu_{{0}} > \mu_{{0}}\mu_{{1}} \frac{-\Lambda_0-\Lambda_1}{\mu_0+\mu_1}+\Lambda_{{0}}\mu_{{1}}+\Lambda_{{1}}\mu_{{0}} = \frac{\Lambda_0 \mu_1^2 + \Lambda_1 \mu_0^2}{\mu_0+\mu_1} >0,
\end{displaymath}

\noindent hence \eqref{neq2} holds. Similarly, assume $\mu_0+\mu_1 < 0$. Then \eqref{neq1} leads to 

\begin{displaymath}
\xi > \frac{-\Lambda_0-\Lambda_1}{\mu_0+\mu_1} > 0.
\end{displaymath}

\noindent Hence, first multiplier in \eqref{neq2} is positive. Consider the second one:

\begin{displaymath}
\xi \mu_{{0}}\mu_{{1}}+\Lambda_{{0}}\mu_{{1}}+\Lambda_{{1}}\mu_{{0}} < \mu_{{0}}\mu_{{1}} \frac{-\Lambda_0-\Lambda_1}{\mu_0+\mu_1}+\Lambda_{{0}}\mu_{{1}}+\Lambda_{{1}}\mu_{{0}} = \frac{\Lambda_0 \mu_1^2 + \Lambda_1 \mu_0^2}{\mu_0+\mu_1} <0,
\end{displaymath}

\noindent hence \eqref{neq2} holds. Obviously, if $\mu_0+\mu_1 = 0$, \eqref{neq1} cannot hold. Proposition is proven.

Inverting Laplace transform in \eqref{Lapl}, we get

\begin{equation}\label{L01}
\begin{split}
& { L_0} \left( t,\xi \right) ={\frac {\Lambda_{{1}}\mu_{{1}}{{\rm e}
^{-\xi\,B}}{{\rm e}^{nt}}\Omega_{{n}} \left( t \right) }{-n+m}}-{
\frac {\Lambda_{{1}}\mu_{{1}}{{\rm e}^{-\xi\,B}}{{\rm e}^{mt}}\Omega_{
{m}} \left( t \right) }{-n+m}}+ \\ & \left( -{\frac {\mu_{{1}}\Psi_n
 \left( t \right) \Lambda_{{1}}}{-n+m}}-{\frac {\mu_{{0}} \left( \xi\,
\mu_{{1}}+n+\Lambda_{{1}} \right) \Phi_{{n}} \left( t \right) }{-n+m}}
 \right) {{\rm e}^{nt}}+ \\ & \left( {\frac {\mu_{{1}}\Psi_m \left( t
 \right) \Lambda_{{1}}}{-n+m}}+{\frac {\mu_{{0}} \left( \xi\,\mu_{{1}}
+m+\Lambda_{{1}} \right) \Phi_{{m}} \left( t \right) }{-n+m}} \right) 
{{\rm e}^{mt}}+ \left( -{\frac {\Lambda_{{1}}{{\rm e}^{nt}}}{-n+m}}+{
\frac {\Lambda_{{1}}{{\rm e}^{mt}}}{-n+m}} \right) { L_1} \left( 0,
\xi \right) + \\& \left( -{\frac { \left( \xi\,\mu_{{1}}+n+\Lambda_{{1}}
 \right) {{\rm e}^{nt}}}{-n+m}}+{\frac { \left( \xi\,\mu_{{1}}+m+
\Lambda_{{1}} \right) {{\rm e}^{mt}}}{-n+m}} \right) { L_0} \left( 0
,\xi \right) ,\\ & { L_1} \left( t,\xi \right) =-{\frac {\Omega_{{m}}
 \left( t \right) \mu_{{1}} \left( \mu_{{0}}\xi+m+\Lambda_{{0}}
 \right) {{\rm e}^{-\xi\,B+mt}}}{-n+m}}+{\frac {\mu_{{1}} \left( \mu_{
{0}}\xi+n+\Lambda_{{0}} \right) {{\rm e}^{-\xi\,B+nt}}\Omega_{{n}}
 \left( t \right) }{-n+m}}+ \\ & \left( -{\frac {\Lambda_{{0}}\mu_{{0}}\Phi
_{{n}} \left( t \right) }{-n+m}}-{\frac {\mu_{{1}} \left( \mu_{{0}}\xi
+n+\Lambda_{{0}} \right) \Psi_n \left( t \right) }{-n+m}} \right) {
{\rm e}^{nt}}+ \\& \left( {\frac {\Lambda_{{0}}\mu_{{0}}\Phi_{{m}} \left( 
t \right) }{-n+m}}+{\frac {\mu_{{1}} \left( \mu_{{0}}\xi+m+\Lambda_{{0
}} \right) \Psi_m \left( t \right) }{-n+m}} \right) {{\rm e}^{mt}}+ \left( {
\frac {{{\rm e}^{mt}}\Lambda_{{0}}}{-n+m}}-{\frac {{{\rm e}^{nt}}
\Lambda_{{0}}}{-n+m}} \right) { L_0} \left( 0,\xi \right)+ \\&
 \left( -{\frac { \left( \mu_{{0}}\xi+n+\Lambda_{{0}} \right) {{\rm e}
^{nt}}}{-n+m}}+{\frac { \left( \mu_{{0}}\xi+m+\Lambda_{{0}} \right) {
{\rm e}^{mt}}}{-n+m}} \right) { L_1} \left( 0,\xi \right),
\end{split}
\end{equation}

\noindent where $\Pi _{{k}} \left( t \right) =\int_{0}^{t}\!\pi  \left( \tau \right) {
{\rm e}^{-k\tau}}\,{\rm d}\tau$ for $\pi \in \{\phi,\psi,\omega\}$ and $k \in \{m,n\}$. Since the process stops with probability one, $L_0$ and $L_1$ tend to zero as $t$ tends to infinity for any $\xi$. In view of Proposition 1, the necessary condition for this is that the coefficients of $\rm e^{m t}$ tend to zero as $t$ tends to infinity. This leads to the system of two equations, one of them turns out to be identity and the second one is

\begin{equation}\label{Bounded}
{ L_1} \left( 0,\xi \right) =-{\frac { \left( \xi\,\mu_{{1}}+m+
\Lambda_{{1}} \right) { L_0} \left( 0,\xi \right) }{\Lambda_{{1}}}}-{\frac {\mu_{{0}} \left( \xi\,\mu_{{1}}+m+\Lambda_{{1}
} \right) \Phi_m}{\Lambda_{{1}}}}  -\mu_{{1}}\Psi_m+\mu_{{1}}{{\rm e}^{-\xi\,B}}\Omega_m,
\end{equation}

\noindent where $\Pi _{{k}} =\int_{0}^{\infty}\!\pi  \left( \tau \right) {
{\rm e}^{-k\tau}}\,{\rm d}\tau$. We now consider the lower boundary condition. Twice integrating \eqref{Ldef} for $s=1$ by parts, we get

\begin{displaymath}
\begin{split}
& 0=\frac{\partial}{\partial A} { F_1} \left( t,0 \right) ={\xi}^{2}{ L_1}
 \left( t,\xi \right) +\xi\,{ F_1} \left( t,B \right) {{\rm e}^{-\xi
\,B}}-\xi\,{ F_1} \left( t,0 \right) +  {{\rm e}^{-\xi\,B}} {\frac {\partial }{
\partial A}}{ F_1} \left( t,B \right) - \\ &
\int_{0}^{B}\! \left( {\frac {\partial ^{2}}{\partial {A}^{2}}}{ F_1
} \left( t,A \right)  \right) {{\rm e}^{-\xi\,A}}\,{\rm d}A.
\end{split}
\end{displaymath}

\noindent Applying mean theorem for the integral, we get

\begin{displaymath}
\begin{split}
& 0= \frac{\partial}{\partial A} { F_1} \left( t,0 \right) ={\xi}^{2}{ L_1}
 \left( t,\xi \right) +\xi\,{ F_1} \left( t,B \right) {{\rm e}^{-\xi
\,B}}-\xi\,{ F_1} \left( t,0 \right) + {{\rm e}^{-\xi\,B}} {\frac {\partial }{
\partial A}}{ F_1} \left( t,B \right)  - \\ & {
\frac {1}{\xi} \frac{\partial^2}{\partial A^2}  { F_1}   \left( t,{
 \hat{A}} \left( t,\xi \right)  \right) }+{\frac{{{\rm e}^{-\xi\,B}}}{{\xi}} { \frac{\partial^2}{\partial A^2} { F_1} \left( t,{ \hat{A}} \left( t,
\xi \right)  \right) }}
\end{split}
\end{displaymath}

\noindent for some point $\hat{A} \left( t,\xi \right)$. Tending $\xi$ to $+\infty$, we get

\begin{equation}\label{Lim}
0=\lim _{\xi\rightarrow \infty } \left({\xi}^{2}{ L_1} \left( t,\xi
 \right) -\xi\,\psi \left( t \right)\right).
\end{equation}

\noindent Substituting \eqref{Bounded} to \eqref{L01}, we obtain

\begin{equation}\label{L1}
\begin{split}
& { L_1} \left( t,\xi \right) =-{\frac {\Omega_{{m}}
 \left( t \right) \mu_{{1}} \left( \mu_{{0}}\xi+m+\Lambda_{{0}}
 \right) {{\rm e}^{-\xi\,B+mt}}}{-n+m}}+{\frac {\mu_{{1}} \left( \mu_{
{0}}\xi+n+\Lambda_{{0}} \right) {{\rm e}^{-\xi\,B+nt}}\Omega_{{n}}
 \left( t \right) }{-n+m}}+ \\ & \left( -{\frac {\Lambda_{{0}}\mu_{{0}}\Phi
_{{n}} \left( t \right) }{-n+m}}-{\frac {\mu_{{1}} \left( \mu_{{0}}\xi
+n+\Lambda_{{0}} \right) \Psi_n \left( t \right) }{-n+m}} \right) {
{\rm e}^{nt}}+ \\& \left( {\frac {\Lambda_{{0}}\mu_{{0}}\Phi_{{m}} \left( 
t \right) }{-n+m}}+{\frac {\mu_{{1}} \left( \mu_{{0}}\xi+m+\Lambda_{{0
}} \right) \Psi_m \left( t \right) }{-n+m}} \right) {{\rm e}^{mt}}+ \left( {
\frac {{{\rm e}^{mt}}\Lambda_{{0}}}{-n+m}}-{\frac {{{\rm e}^{nt}}
\Lambda_{{0}}}{-n+m}} \right) { L_0} \left( 0,\xi \right)+ \\&
 \left( -{\frac { \left( \mu_{{0}}\xi+n+\Lambda_{{0}} \right) {{\rm e}
^{nt}}}{-n+m}}+{\frac { \left( \mu_{{0}}\xi+m+\Lambda_{{0}} \right) {
{\rm e}^{mt}}}{-n+m}} \right) \\& \left( -{\frac { \left( \xi\,\mu_{{1}}+m+
\Lambda_{{1}} \right) { L_0} \left( 0,\xi \right) }{\Lambda_{{1}}}}-{\frac {\mu_{{0}} \left( \xi\,\mu_{{1}}+m+\Lambda_{{1}
} \right) \Phi_m}{\Lambda_{{1}}}} -\mu_{{1}}\Psi_m+\mu_{{1}}{{\rm e}^{-\xi\,B}}\Omega_m \right).
\end{split}
\end{equation}

\noindent Now we apply the following representation of every integral:

\begin{displaymath}
\int_{0}^{T}\!\pi  \left( \tau \right) {{\rm e}^{-k\tau}}\,{\rm d}\tau
=-{\frac {\pi  \left( T \right) {{\rm e}^{-kT}}}{k}}+{\frac {\pi 
 \left( 0 \right) }{k}}-{\frac {\pi'(T) {{\rm e}^{-kT}}}{{k}^{2}}}+{\frac {\pi'  \left( 0 \right) }{{k}^{2}}}+{\frac { \pi''  \left( \theta \right) }{{k}^{
3}}}-{\frac { \pi'' \left( \theta \right) {{\rm e}^{-kT}}}{{k}^{3}}}, T \in \{t,\infty\}.
\end{displaymath}

Substituting it to \eqref{L1} and applying Taylor series for big $\xi$, we get

\begin{equation}\label{L1appr}
{ L_1} \left( t,\xi \right) \approx {\frac {\psi \left( t \right) \xi\,\mu_
{{1}}-\psi \left( t \right) \Lambda_{{1}}+\phi \left( t \right) 
\Lambda_{{0}}-\psi' \left( t \right) }{\mu_{{
1}}{\xi}^{2}}}
\end{equation}

\noindent Substituting \eqref{L1appr} to \eqref{Lim} we get

\begin{equation}\label{Diff}
\psi \left( t \right) \Lambda_{{1}}-\phi \left( t \right) \Lambda_{{0}
}+\psi' \left( t \right) = 0.
\end{equation}

Substituting \eqref{Diff} to \eqref{Bounded} and integrating by parts, we get

\begin{equation}\label{L11}
\begin{split}
& { L_1} \left( 0,\xi \right) = \\ & {\frac {\mu_{{0}}\mu_{{1}}\xi\,\psi
 \left( 0 \right) -{ L_0} \left( 0,\xi \right) \xi\,\Lambda_{{0}}\mu
_{{1}}+\mu_{{0}}m\psi \left( 0 \right) +\mu_{{0}}\psi \left( 0
 \right) \Lambda_{{1}}-{ L_0} \left( 0,\xi \right) m\Lambda_{{0}}-{
 L_0} \left( 0,\xi \right) \Lambda_{{1}}\Lambda_{{0}}}{\Lambda_{{0}}
\Lambda_{{1}}}}- \\ & {\frac { \left( m\xi\,\mu_{{0}}\mu_{{1}}+\xi\,\Lambda_
{{1}}\mu_{{0}}\mu_{{1}}+{m}^{2}\mu_{{0}}+2\,m\Lambda_{{1}}\mu_{{0}}+
\Lambda_{{0}}\Lambda_{{1}}\mu_{{1}}+\Lambda_1^{2}\mu_{{0}}
 \right) \int_{0}^{\infty }\!\psi \left( \tau \right) {{\rm e}^{-m\tau
}}\,{\rm d}\tau}{\Lambda_{{0}}\Lambda_{{1}}}}+ \\ & \mu_{{1}}{{\rm e}^{-\xi
\,B}}\int_{0}^{\infty }\!\omega \left( \tau \right) {{\rm e}^{-m\tau}}
\,{\rm d}\tau.
\end{split}
\end{equation}

We now exclude $\xi$ from \eqref{L11}. To do so, we express $\xi$ through $m$ from \eqref{mn}. It can be done in two ways:

\begin{equation}\label{xi}
\begin{split}
& \xi_1 = {\frac {-\mu_{{0}}m-m\mu_{{1}}-\Lambda_{{0}}\mu_{{1}}-\Lambda_{{
1}}\mu_{{0}}+\sqrt{r}}{2 \mu_{{0}}\mu_{{1}}}},\xi_2 = -{\frac {\mu_{{0}}m+m\mu_{{1
}}+\Lambda_{{0}}\mu_{{1}}+\Lambda_{{1}}\mu_{{0}}+\sqrt{r}}{2\mu_{{0}}\mu_{{1}}
}}, \\ & r={m}^{2}\mu_0^{2}-2\,{m}^{2}\mu_{{0}}\mu_{{1}}+{m}^{2
}\mu_1^{2}-2\,m\Lambda_{{0}}\mu_{{0}}\mu_{{1}}+2\,m\Lambda_{{0}}
\mu_1^{2}+2\,m\Lambda_{{1}}\mu_0^{2}-2\,m\Lambda_{{1}}\mu_
{{0}}\mu_{{1}}+ \\ & \Lambda_0^{2}\mu_1^{2}+2\,\Lambda_{{0}}
\Lambda_{{1}}\mu_{{0}}\mu_{{1}}+\Lambda_1^{2}\mu_0^{2}.
\end{split}
\end{equation}

\noindent We also introduce

\begin{equation}\label{UV}
U=m\mu_{{0}}+m\mu_{{1}}+\Lambda_{{0}}\mu_{{1}}+\Lambda_{{1}}\mu_{{0}}
-\sqrt {r},W=m\mu_{{0}}+m\mu_{{1}}+\Lambda_{{0}}\mu_{{1}}+\Lambda_{{1}
}\mu_{{0}}+\sqrt {r}.
\end{equation}

\noindent Substituting \eqref{xi} into \eqref{L11}, we get:

\begin{equation}\label{Sys}
\begin{split}
& { L_1} \left( 0,-{\frac {U}{2\mu_{{0}}\mu_{{1}}}} \right) =-
\,{\frac {-2\,m\mu_{{0}}\psi \left( 0 \right) -2\,\mu_{{0}}\psi
 \left( 0 \right) \Lambda_{{1}}+U\psi \left( 0 \right) }{2\Lambda_{{0}}
\Lambda_{{1}}}}+ \\ & \mu_{{1}}{{\rm e}^{{\frac {UB}{2\mu_{{0}}\mu_{{1}}
}}}}\int_{0}^{\infty }\!\omega \left( \tau \right) {{\rm e}^{-m\tau}}
\,{\rm d}\tau+{\frac {-2\,m\mu_{{0}}-2\,\Lambda_{{1}}\mu_{{0}}+U
}{2\Lambda_{{1}}\mu_{{0}}}{ L_0} \left( 0,-{\frac {U}{2\mu_{{0}}
\mu_{{1}}}} \right) }+ \\ & {\frac { \left( -2\,{m}^{2}\mu_{{0}}-4\,m
\Lambda_{{1}}\mu_{{0}}-2\,\Lambda_{{0}}\Lambda_{{1}}\mu_{{1}}-2 
\Lambda_1^{2}\mu_{{0}}+Um+U\Lambda_{{1}} \right) \int_{0}^{
\infty }\!\psi \left( \tau \right) {{\rm e}^{-m\tau}}\,{\rm d}\tau}{
2\Lambda_{{0}}\Lambda_{{1}}}}, \\ &
{ L_1} \left( 0,-{\frac {W}{2 \mu_{{0}}\mu_{{1}}}} \right) =-
\,{\frac {-2\,m\mu_{{0}}\psi \left( 0 \right) -2\,\mu_{{0}}\psi
 \left( 0 \right) \Lambda_{{1}}+W\psi \left( 0 \right) }{2\Lambda_{{0}}
\Lambda_{{1}}}}+ \\ & \mu_{{1}}{{\rm e}^{{\frac {W B}{2\mu_{{0}}\mu_{{1}}
}}}}\int_{0}^{\infty }\!\omega \left( \tau \right) {{\rm e}^{-m\tau}}
\,{\rm d}\tau+{\frac {-2\,m\mu_{{0}}-2\,\Lambda_{{1}}\mu_{{0}}+W
}{2\Lambda_{{1}}\mu_{{0}}}{ L_0} \left( 0,-{\frac {W}{2\mu_{{0}}
\mu_{{1}}}} \right) }+ \\ & {\frac { \left( -2\,{m}^{2}\mu_{{0}}-4\,m
\Lambda_{{1}}\mu_{{0}}-2\,\Lambda_{{0}}\Lambda_{{1}}\mu_{{1}}-2 
\Lambda_1^{2}\mu_{{0}}+W m+W\Lambda_{{1}} \right) \int_{0}^{
\infty }\!\psi \left( \tau \right) {{\rm e}^{-m\tau}}\,{\rm d}\tau}{
2\Lambda_{{0}}\Lambda_{{1}}}}.
\end{split}
\end{equation}

This is the linear system of equations on two unknowns. We derive its solution as the combination of fundamental solutions. Indeed, let $\psi_A^0 (t)$ and $\omega_A^0 (t)$ be solutions of \eqref{Sys} for $F_0 (0,A) = 0$, $F_1 (0,A) = \delta_A (x)$ and $\psi(0) = 0$ and $\psi_A^1 (t)$ and $\omega_A^1 (t)$ be solutions of \eqref{Sys} for $F_0 (0,A) = \delta_A (x)$, $F_1 (0,A) = 0$ and $\psi(0) = 0$ and $\psi_A^2 (t)$ and $\omega_A^2 (t)$ be solutions of \eqref{Sys} for $F_0 (0,A) = 0$, $F_1 (0,A) = 0$. Then the general solution of \eqref{Sys} can be found as

\begin{equation}\label{psiomega}
\begin{split}
& \psi(t) = \int_0^B \psi_A^0 (\tau) F_0 (0,A) d A + \int_0^B \psi_A^1 (\tau) F_1 (0,A) d A + \psi^2 (\tau), \\ & \omega(t) = \int_0^B \omega_A^0 (\tau) F_0 (0,A) d A + \int_0^B \omega_A^1 (\tau) F_1 (0,A) d A + \omega^2 (\tau).
\end{split}
\end{equation}

Substituting $F_0 (0,A) = 0$, $F_1 (0,A) = \delta_A (x)$ and $\psi(0) = 0$ into \eqref{Sys} and solving the system, we get

\begin{equation}\label{fin0}
\begin{split}
& \int_{0}^{\infty }\!\omega_A^0 \left( \tau \right) {{\rm e}^{-m\tau}}
\,{\rm d}\tau=\frac{ \left( {{\rm e}^{{\frac {U A}{2\mu_
{{0}}\mu_{{1}}}}}}U\Lambda_{{0}}-{{\rm e}^{{\frac {W A
}{2\mu_{{0}}\mu_{{1}}}}}}W\Lambda_{{0}}
 \right)}{\mu_0  \left({{\rm e}^{{\frac {W B}{2\mu_{{0}}\mu_{{1}}}}}}S_{{
2}}+S_{{1}}{{\rm e}^{{\frac {U B}{2\mu_{{0}}\mu_{{1}}}}}} \right)},
 \\ & 
 \int_{0}^{\infty }\!\psi_A^0 \left( \tau \right) {{\rm e}^{-m\tau}}
\,{\rm d}\tau={\frac {2\Lambda_{{1}}\Lambda_{{0}}\mu_{{1}}\int_{0}^{
\infty }\!\omega^A_0 \left( \tau \right) {{\rm e}^{-m\tau}}\,{\rm d}\tau}{S
_{{2}}}{{\rm e}^{{\frac {UB}{2\mu_{{0}}\mu_{{1}}}}}}}+ \\ & -{\frac {\Lambda_{{0
}} \left( m\mu_{{0}}-m\mu_{{1}}-\Lambda_{{0}}\mu_{{1}}+\Lambda_{{1}}
\mu_{{0}}+s \right) }{S_{{2}}\mu_{{0}}}{{\rm e}^{{\frac {UA}{2\mu_
{{0}}\mu_{{1}}}}}}}.
\end{split}
\end{equation}

Substituting $F_0 (0,A) = \delta_A (x)$, $F_1 (0,A) = 0$ and $\psi(0) = 0$ into \eqref{Sys} and solving the system, we get

\begin{equation}\label{fin1}
\begin{split}
& \int_{0}^{\infty }\!\omega_A^1 \left( \tau \right) {{\rm e}^{-m\tau}}
\,{\rm d}\tau={ \frac{\left({{\rm e}^{{\frac {W A}{2\mu_{{0}}\mu_{{1}}}}}}S_{{2}}+{
{\rm e}^{{\frac {U A}{2 \mu_{{0}}\mu_{{1}}}}}}S_{{1}} \right)}{ 
\mu_{{1}} \left( {{\rm e}^{{\frac {W B}{2 \mu_{{0}}\mu_{{1}}}}}}S_{{2}}+S_{{1
}}{{\rm e}^{{\frac {U B}{2 \mu_{{0}}\mu_{{1}}}}}} \right)} }, \\ & 
\int_{0}^{\infty }\!\psi^1_A \left( \tau \right) {{\rm e}^{-m\tau}}
\,{\rm d}\tau={\frac {2\Lambda_{{1}}\Lambda_{{0}}\mu_{{1}}\int_{0}^{
\infty }\!\omega_A^1 \left( \tau \right) {{\rm e}^{-m\tau}}\,{\rm d}\tau}{S
_{{2}}}{{\rm e}^{{\frac {U B}{2 \mu_{{0}}\mu_{{1}}}}}}}-{\frac {
2 \Lambda_{{1}}\Lambda_{{0}}}{S_{{2}}}{{\rm e}^{{\frac {U A}{2 \mu_{{0
}}\mu_{{1}}}}}}}.
\end{split}
\end{equation}

Substituting $F_0 (0,A) = 0$ and $F_1 (0,A) = 0$ into \eqref{Sys} and solving the system, we get

\begin{equation}\label{fin2}
\begin{split}
& \int_{0}^{\infty }\!\omega^2 \left( \tau \right) {{\rm e}^{-m\tau}}
\,{\rm d}\tau = \frac{2 s \psi(0)}{{ 
 \left( {{\rm e}^{{\frac {W B}{2 \mu_{{0}}\mu_{{1}}}}}}S_{{2}}+S_{{1
}}{{\rm e}^{{\frac {U B}{2 \mu_{{0}}\mu_{{1}}}}}} \right)}},  
\\ & \int_{0}^{\infty }\!\psi^2 \left( \tau \right) {{\rm e}^{-m\tau}}
\,{\rm d}\tau = \\ & \frac{-\psi \left( 0 \right)  \left( -2\,m\mu_{{0}}-2\,\Lambda_{{1}}\mu_{{0}
}+U \right) {{\rm e}^{{\frac {W B}{2 \mu_{{0}}\mu_{{1}}}}}}+\psi
 \left( 0 \right)  \left( -2\,m\mu_{{0}}-2\,\Lambda_{{1}}\mu_{{0}}+W
 \right) {{\rm e}^{{\frac {U B}{2 \mu_{{0}}\mu_{{1}}}}}}}{{  \left( {{\rm e}^{{\frac {W B}{2 \mu_{{0}}\mu_{{1}}}}}}S_{{2}}+S_{{1
}}{{\rm e}^{{\frac {U B}{2 \mu_{{0}}\mu_{{1}}}}}} \right)}}.
\end{split}
\end{equation}

\noindent Hence, in order to derive explicit formulae for $\omega_A^0, \psi_A^0, \omega_A^1, \psi_A^1,\omega^2, \psi^2$, we need to find corresponding inverse Laplace transforms in \eqref{fin0}, \eqref{fin1} and \eqref{fin2}. After that $\omega$ and $\psi$ can be found from \eqref{psiomega} and $\phi$ can be found from \eqref{Diff}. After that $L_0$ and $L_1$ can be found from \eqref{L01}, which gives the full description of the dynamics of the process.

\bibliographystyle{abbrv}
\bibliography{MyBibliography}

\end{document}